\documentclass[12pt]{amsart}  
\usepackage{amssymb}
\usepackage{mathptmx}

\newcommand{\mm}{\mathfrak m}  
\newcommand{\Z}{\mathbb{Z}}  
  
\newcommand{\N}{\mathbb{N}}

\newcommand {\PP}{\mathbb{P}}
\DeclareMathOperator{\depth}{depth} 
\DeclareMathOperator{\GL}{GL}  
\DeclareMathOperator{\Ext}{Ext}

\DeclareMathOperator{\gin}{gin}

\DeclareMathOperator{\ini}{in}

\newcommand{\fp}{{\mathfrak p}}

\DeclareMathOperator{\dirsum}{\oplus}

\DeclareMathOperator{\Dirsum}{\bigoplus}  
  
\DeclareMathOperator{\pnt}{\raise 0.5mm \hbox{\large\bf.}}

\DeclareMathOperator{\Deg}{Deg}  
\DeclareMathOperator{\adeg}{adeg}  
\DeclareMathOperator{\hdeg}{hdeg}  
\DeclareMathOperator{\sdeg}{sdeg}  
\newcommand{\Ex}[3]{{\rm Ext}_{#2}^{#1}(#3,\omega_{#2})}
\newtheorem{thm}{\bf Theorem}[section]   
\newtheorem{lem}[thm]{\bf Lemma}  
\newtheorem{cor}[thm]{\bf Corollary}  
\newtheorem{prop}[thm]{\bf Proposition}

\theoremstyle{definition}
\newtheorem{defn}[thm]{\bf Definition}  
  
\newtheorem{rem}[thm]{\bf Remark}  
\newtheorem{ex}[thm]{\bf Example}   
\newtheorem{alg}[thm]{\bf Algorithm}   
\title{Extended degree functions and monomial modules}  
  
\author[Uwe Nagel]{Uwe Nagel$^*$} 
\address{Department of Mathematics, University of Kentucky, Lexington, KY 40506-0027, USA}
\email{uwenagel@ms.uky.edu}

\author{Tim R\"omer} 
\address{FB Mathematik/Informatik, Universit\"at Osnabr\"uck, 49069 Osnabr\"uck,  Germany}
\email{troemer@mathematik.uni-osnabrueck.de} 

\thanks{$^*$ The first author gratefully acknowledges partial support by 
a Special Faculty Research Fellowship from the University of Kentucky.}

\begin{document}

\begin{abstract}  
The arithmetic degree, the smallest extended degree, and the homological degree are invariants that have been proposed as alternatives 
of the degree of a module if this module is not Cohen-Macaulay. We compare these degree functions and study their behavior when passing 
to the generic initial or the lexicographic submodule. This leads to various bounds and to counterexamples to a conjecture of Gunston and 
Vasconcelos, respectively. Particular attention is given to the class of sequentially Cohen-Macaulay modules. The results in this case lead to an 
algorithm that computes the smallest extended degree. 
\end{abstract}
  
\maketitle 

\tableofcontents

%
%
%
\section{Introduction}  

Let $M$ be a finitely generated graded module over the polynomial ring $S$. If $M$ is Cohen-Macaulay then several invariants of $M$ can be bounded 
using the degree of $M$. This is no longer true if $M$ is not Cohen-Macaulay. In this case, one tries to replace the degree of $M$ by an invariant 
that better captures the structure of $M$. One such invariant is the {\em arithmetic degree} (cf.\ \cite{BAMU}) of $M$  
$$
\adeg M = \sum \deg M_{\fp} \cdot  \deg \fp
$$  
where the sum runs over the associated prime ideals of $M$. 

More recently, Vasconcelos \cite{VA98} has axiomatically introduced so-called extended degrees (cf.\ Section \ref{sec-functions}). 
They are designed to provide measures for the size and the complexity of the structure of $M$. 
The first concrete example of an extended degree is Vasconcelos' {\em homological degree} \cite{VA98}. It is recursively defined by 
$$
\hdeg M = \deg M  + \sum_{i=0}^{d-1} \binom{d-1}{i} 
\hdeg \Ext_S^{n-i}(M,S) 
$$ 
where $d := \dim M$. 
Gunston (\cite{GU}, cf.\ also \cite{NA03}, Lemma 4.2) has shown that among all extended degrees there is a minimal one which we just 
call the {\em smallest extended degree} $\sdeg M$. In this paper we compare these three degrees and study their behavior when we replace $M$ 
by a related monomial module. This leads to various bounds. 

A difficulty when dealing with the smallest extended degree is that, in general, there is no formula to compute it. 
However, we show that such a formula does exist if $M$ is either a sequentially Cohen-Macaulay (cf.\ Section \ref{sec-seqCM}) 
or a Buchsbaum module (cf.\ Section \ref{sec_buchsbaum}). As a first application of these formulas, we show in Section \ref{sec-bounds}  
that every module $M$ satisfies
$$
\deg M \leq \adeg M \leq \sdeg M \leq \hdeg M. 
$$
This refines Vasconcelos' Proposition 9.4.2 in \cite{VA}.
Moreover, our formulas show that 
$$
\adeg M = \sdeg M
$$
if $M$ is sequentially Cohen-Macaulay, and that 
$$
\sdeg M = \hdeg M 
$$
if $M$ is a Buchsbaum module. 

The case of sequentially Cohen-Macaulay modules is of particular importance because such modules naturally occur. 
Indeed, write $M = F/U$ where $F$ is a free $S$-module and $U \subset F$ is a graded submodule. By now it is a standard technique to draw 
conclusions about $F/U$ by considering $F/\gin (U)$ where $\gin(U)$ is the generic initial module of $U$ with respect to the reverse lexicographic 
order on $F$ (cf.\ Section \ref{sec-functions}). In order to get bounds for invariants on $M$ that depend on its Hilbert function, it is often 
useful to compare $M = F/U$ with $F/U^{lex}$ where $U^{lex}$ is the lexicographical submodule of $F$ that has the same Hilbert 
function as $U$ (cf.\ Section \ref{sec-functions}). Both, $\gin(U)$ and $U^{lex}$ are 
{\em Borel-fixed} (cf.\ Section \ref{sec-functions}), 
thus they are sequentially Cohen-Macaulay (cf.\ Lemma \ref{lem-borel-t-is-seq-CM}). 
Hence, our formulas apply and we use them to show that we have for every module $M = F/U$ 
$$
\adeg F/U \leq \adeg F/\gin(U) \leq \adeg F/U^{lex}
$$
and 
$$
\sdeg F/U = \sdeg F/\gin(U) \leq \sdeg F/U^{lex}. 
$$
Note that the first inequality for $\adeg$ extends a result of Sturmfels, Trung, and Vogel \cite[Theorem 2.3]{STTRVO} 
from ideals to submodules whereas the equality for $\sdeg$ is due to \cite{GU} (cf.\ also \cite{NA03}). 
In spite of the estimates above, it is natural to conjecture (cf.\ \cite{GU} and \cite[page 262]{VA}) 
that  we have for every module $M = F/U$ either  always the relation 
$$
\hdeg F/U \geq \hdeg F/\gin (U) 
$$
or 
$$
\hdeg F/U \leq \hdeg F/\gin (U).
$$
Since it is often possible to compare invariants of $F/U$ and $F/U^{lex}$, one might also suspect that there is always either the relation  
$$
\hdeg F/U \geq \hdeg F/U^{lex}
$$
or 
$$
\hdeg F/U \leq \hdeg F/U^{lex}.
$$
In fact, this work began as an attempt to prove these conjectures. Somewhat surprisingly we show in Section \ref{sec-examples} 
that none of the conjectured relations is always true by exhibiting suitable modules. 

Our formulas and estimates for the degree functions are in terms of the degrees of certain extension modules. In the final section, 
we show that these degrees can very efficiently be computed in case of monomial modules of Borel-type. As a consequence, we get 
a fast algorithm for computing $\sdeg F/U$ provided we know $\gin (U)$. 

Throughout the paper we consider finitely generated graded modules over the polynomial ring $S$. 
However, using \cite{CN} our results for sequentially Cohen-Macaulay and Buchsbaum modules remain valid for modules over an arbitrary 
Noetherian local ring $(R, \mm)$ provided monomial modules are not involved. In the latter case, the result are still true for modules 
over a regular local ring $(R, \mm)$ of dimension $n$ where the maximal ideal $\mm$ is generated by $x_1,\ldots,x_n$.

%
%
%
\section{Degree functions}  \label{sec-functions} 

In this section we introduce several degree functions of modules.
We briefly recall definitions and notation used in this paper. 
For unexplained terminology we refer to the book of Bruns and Herzog \cite{BRHE98}.

Throughout this paper $K$ is always an infinite field and
$S=K[x_1,\dots,x_n]$ is the polynomial ring over $K$ with its standard grading where 
 $\deg x_i=1$ for $i=1,\dots,n$.
We denote by $\mm=(x_1,\dots,x_n)$ the unique graded maximal ideal of $S$. 
A standard graded $K$-algebra $R$ is of the form $S/I$
for a graded ideal $I \subset S$. 

Usually we denote by $M$ a finitely generated graded $S$-module
of dimension $d=\dim M$. Its $i$-th local cohomology module is denoted by $H^i_{\mm} (M)$. 
The {\em Hilbert function} $H_M$ of $M$ is defined by
$$
H_M \colon \Z \to \N, j\mapsto \dim_K M_j.
$$
It is well-known that 
there exists a polynomial $P_M$ of degree $d -1$
such that for $j\gg 0$ we have that $H_M(j)=P_M(j)$.
We write 
$$
P_M(t)= \frac{e(M)}{(d-1)!} t^{d-1} + \dots \text{ (terms of lower degree)}.
$$
We define the {\em degree} $\deg M$ (or {\em multiplicity} of $M$) to be $e(M)$ 
if $d>0$
and $\deg M = l(M)$ if $d=0$. Here $l(\cdot)$ denotes the length function of a 
module $M$. 
The degree of $M$ has many nice properties, especially if $M$ is
a Cohen-Macaulay module (CM module for short).

There are several attempts do
define degree functions for a module $M$
that coincide with the  degree if $M$ is a CM module,
but also have nice properties for non-CM modules.
We refer to the nice book of Vasconcelos \cite{VA98} for details on this subject. 
One such proposal is due to Bayer and Mumford who introduced in \cite{BAMU} the 
{\em arithmetic degree} that has been studied by several authors  in the last decades  (see, e.g., \cite{MIVO}, 
\cite{STTRVO}, or \cite{VA}). Vasconcelos \cite[Proposition 9.1.2]{VA} has shown that the arithmetic degree 
can be computed using the formula 
$$
\adeg M= \sum_{i=0}^n  \deg \Ext_S^{i}({\Ext_S^{i}(M,\omega_S)},\omega_S) 
$$  
where $\omega_S=S(-n)$ is the canonical module of $S$. 
But there are some disadvantages.
For example, if $y \in S_1$ is an {\em $M$-regular element}, i.e.\ it is a non-zero 
divisor of $M$ 
(sometimes also called a regular hyperplane section),
then
$$
\adeg M \leq \adeg M/yM.
$$
But if a degree function reflects the complexity of the module,
then $M/yM$ should have a smaller degree than $M$.

In \cite{VA98}, Vasconcelos axiomatically 
defined the following concept. A
numerical function 
$\Deg$ that assigns to every finitely generated graded $S$-module  a non-negative 
integer is said to be an {\em extended degree function} if
it satisfies the following conditions:
\begin{enumerate}
\item
If $L=H^0_{\mm}(M)$, then $\Deg M =\Deg M/L +l(L).$
\item
If $y\in S_1$ is sufficiently general and $M$-regular, 
then $\Deg M \geq \Deg M/yM.$
\item
If $M$ is a CM module, then $\Deg M=\deg M.$
\end{enumerate}

The first example of such an extended degree function has been introduced by Vasconcelos.
The {\em homological degree} of $M$ 
is defined recursively as 
$$
\hdeg M =\deg M  + \sum_{i=0}^{d-1} \binom{d-1}{i} \hdeg \Ext_S^{n-
i}(M,\omega_S).
$$
Note that this is well-defined because  
$\dim \Ext_S^{n-i}(M,\omega_S) \leq i$ for $i=0,\dots,n$.
In \cite{VA98} it is shown that $\hdeg M $ is indeed an extended degree
function. 

Another extended degree function was defined by Gunston in his thesis \cite[Theorem 3.1.2]{GU}, 
the {\em smallest extended degree} $\sdeg$.
Let us recall its axiomatic description.
\begin{thm}
\label{sdeg_define}
There is a unique numerical function $\sdeg$ defined on finitely generated 
graded $S$-modules, satisfying the following conditions:
\begin{enumerate}
\item
If $L=H^0_{\mm}(M)$ then $\sdeg M =\sdeg M/L +l(L)$.
\item
If $y\in S_1$ is sufficiently general 
and $M$-regular, then $\sdeg M =\sdeg M/yM $.
\item
$\sdeg(0)=0$.
\end{enumerate}
\end{thm}

We recall important properties of the function $\sdeg$. (See \cite{NA03} for 
details.)
\begin{enumerate}
\item
$\sdeg$ is indeed an extended degree function.  
For any other extended degree function
$\Deg$ we have that $\sdeg M \leq \Deg M $ for all finitely generated graded 
$S$-modules $M$.
\item
Let $F$ be a finitely generated graded free $S$-module and $U \subset F$ a 
graded submodule. 
Then $\sdeg F/U =\sdeg F/\gin(U)$ where $\gin(U)$ is the {\em generic initial 
submodule} of $U$
with respect to the reverse lexicographic term order on $F$. 
\end{enumerate}

We briefly recall the construction of $\gin(U)$ because we need this module
several times in this paper. For details see, for example, Eisenbud's book  
\cite{EI}.
Let $e_1,\dots,e_m$ be a homogeneous basis for the free graded $S$-module $F$.
For a monomial $x^u \in S$ 
we call an element $x^ue_j$ a {\em monomial} in $F$. The (degree) reverse 
lexicographic term-order $<$ (revlex order for short)
is defined as follows:
$x^ue_s< x^ve_t$ 
if either $\deg x^ue_s< \deg x^ve_t$
or 
$\deg x^ue_s = \deg x^ve_t$ and $x^u < x^v$ in the usual revlex term-order on 
$S$
or
$\deg x^ue_s = \deg x^ve_t$, $x^u =x^v$ and $s>t$.

Consider 
$\GL(n)$ as the group of $K$-linear graded automorphisms of $S$
and let 
$\GL(F)$ be the group of $S$-linear graded automorphisms of $F$.
Then $\mathcal{G}=\GL(n)\ltimes \GL(F)$ acts on $F$ through 
$K$-linear graded automorphisms.
Recall that a {\em monomial submodule} of $F$  
is a module generated by monomials of $F$.
There exists a non-empty open
set $\mathcal{U} \subset \mathcal{G}$ and a unique monomial submodule $U'\subset 
F$ 
with $U'=\ini_>(g(U))$ for every $g \in \mathcal{U}$ 
with respect to the revlex order.
We call $U'$ the generic initial module of $U$ and denote it by $\gin(U)$. 

We will also consider the {\em lexicographic submodule} $U^{lex}$ associated to $U \subset F$. 
The lexicographic order on $F$ is defined by $x^u e_i > x^v e_j$ if either $i < j$ or $i = j$ and the first 
non-zero entry of $u - v$ is positive. A lexicographic submodule is a monomial submodule $V \subset F$ such 
that, for every $i$, $V_i$ is spanned by the first $\dim_K V_i$ monomials of $F_i$ in the lexicographic order. 
If $U \subset F$ is any graded submodule then $U^{lex}$ is the lexicographic submodule of $F$ such that 
$\dim_K U_i = \dim_K (U^{lex})_i$ for all integers $i$. 

Later on we will use the fact that $U^{lex}$ and $\gin (U)$ are Borel-fixed submodules (cf.\ \cite{EI}). 

%
%

%
%
%
\section{Sequentially Cohen-Macaulay modules} 
\label{sec-seqCM}

In this section 
we derive formulas for degree functions when they are restricted to the class of    
{\em sequentially Cohen-Macaulay modules}. 
 The methods developed here
will be very useful in later sections.

Let us briefly recall the definition and some facts about sequentially Cohen-Macaulay 
modules. 
Let $K$ be field and let $R$ be a standard graded Cohen-Macaulay $K$-algebra of 
dimension $n$ with canonical module $\omega_R$.  
The following definition is due to Stanley \cite{ST96}.
\begin{defn}
Let $M$ be a finitely generated graded $R$-module. 
The module $M$ is said to be 
{\em sequentially Cohen-Macaulay} 
(sequentially CM modules for short),
if there exists a finite filtration
\begin{equation}
\label{filter}
0=M_0 \subset M_1 \subset M_2 \subset \dots \subset M_r=M
\end{equation}
of $M$ by graded submodules of $M$ 
such that each quotient $M_i/M_{i-1}$ is Cohen-Macaulay and 
$\dim M_1/M_0 <\dim M_2/M_1 <\dots<\dim M_r/M_{r-1}$.
\end{defn}

We recall some results from \cite{HESB} and \cite{HEPOVL}: 
\begin{enumerate}
\item
The filtration (\ref{filter}) of a sequentially CM module is uniquely determined
and is called the {\em CM-filtration} of $M$.
\item
Setting $d_i= \dim M_i/M_{i-1}$ we have 
$d_i=\dim M_i$ for $i=1,\dots,r$. Furthermore $\dim M=d_r$ and $\depth M=d_1$.
\item
$M$ is sequentially CM if and only if for all $i=0,\ldots, \dim M $ 
we have that the modules $\Ext_R^{n-i}(M,\omega_R)$ are either $0$ or CM of dimension $i$.
In this case, 
\begin{eqnarray*}
\Ext_R^{n-d_i}(M,\omega_R) &\cong& \Ext_R^{n-d_i}(M_i/M_{i-1},\omega_R) 
\text{ for } i=1,\dots,r \text{ and } \\
\Ext_R^{n-i}(M,\omega_R)&=&0 \text{ for } i \neq \{d_1,\dots,d_r\}.
\end{eqnarray*}
\item
A finite direct sum of sequentially CM  modules 
is sequentially CM.
\item
$M$ is sequentially CM if and only if $M/H^0_{\mm}(M)$
is sequentially CM.
\item
Let $y \in R_1$ be an
$M$-regular element that is also regular on all
$\Ext_R^{i}(M,\omega_R)$. Then $M$ is sequentially CM if and only if $M/yM$ is 
sequentially CM.
\end{enumerate}


From now on all modules are assumed to be finitely generated graded modules over $S=K[x_1,\ldots,x_n]$. 

Sequentially Cohen-Macaulay modules occur frequently. We set 
$$
U : I^{\infty} := \bigcup_{k \geq 0} U :_F I^k
$$
if $U$ is a submodule of the free $S$-module $F$ and $I \subset S$ is an ideal. Then we have:

\begin{rem}
\label{boreltype}
Let $I \subset S$ be a graded ideal. Recall the following definition
from \cite{HEPOVL}. The ideal 
$I$ is said to be of {\em Borel type},
if  we have for $i=1,\dots,n$ that
$$
I: x_i^{\infty}= I:(x_1,\dots,x_i)^\infty.
$$
A {\em Borel-fixed} ideal is of Borel-type (see \cite[Proposition 15.24]{EI}), 
hence so is the generic initial ideal $\gin(I)$ 
with respect to the reverse lexicographic term order 
of $I$. In \cite{HEPOVL} Herzog, Popescu and Vladoiu proved that if $I$ is a 
monomial ideal of Borel-type,
then $R=S/I$ is sequentially CM. 
\end{rem} 

Observe that the last result is no longer true if $I$ is not a monomial ideal as the following example shows. 

\begin{ex}  \label{ex-mon-nec} Consider the ideal 
$$
I=(x_1^2,x_1x_2,x_2^2,x_1x_3+x_2x_4) \subset K[x_1,x_2,x_3,x_4].
$$
It defines a double line in $\PP^3$. Let $R=K[x_1,x_2,x_3,x_4]/I$. 
Then we have  (cf., e.g., \cite{NNS}) that
$\dim R=2$ and 
$\Ext^{n-1}_S(R,\omega_S) \cong K$.
Hence $R$ is not sequentially CM, but it is of Borel type because 
\begin{eqnarray*}
I:x_4^\infty &=& I: (x_1,x_2,x_3,x_4)^\infty = I,\\
I:x_3^\infty &=& I: (x_1,x_2,x_3)^\infty = I,\\
I:x_2^\infty &=& I: (x_1,x_2)^\infty = K[x_1,x_2,x_3,x_4],\\
I:x_1^\infty &=& K[x_1,x_2,x_3,x_4].
\end{eqnarray*} 
\end{ex}

The notion of monomial ideals of Borel-type can easily be generalized to 
modules.
Let $F$ be a finitely generated free graded $S$-module 
with homogeneous basis $e_1,\dots,e_m$ and let $U \subseteq F$ be a graded 
submodule.
The module $U$ is said to be of {\em Borel-type} if
$$
U: x_i^{\infty}= U:(x_1,\dots,x_i)^\infty \text{ for } i=1,\dots,n.
$$

As for ideals, we have: 

\begin{lem} \label{lem-borel-t-is-seq-CM} 
If $U \subset F$ is monomial and of Borel-type then $F/U$ is sequentially CM. 

In particular, $F/U$ is sequentially CM if $U = \gin(V)$ or $U = W^{lex}$ for graded submodules $V, W  \subset F$. 
\end{lem}  

\begin{proof} 
By assumption, we can write $U=\dirsum_{i=1}^m I_j e_j $
for monomial ideals $I_j \subset S$ of Borel-type. 
Since, by \cite{HEPOVL}, $S/I_j$ is sequentially CM for $j=1,\dots,m$, 
we have that
$F/U\cong \dirsum_{j=1}^m S/I_j$ is sequentially CM 
because a direct sum of sequentially CM modules is
sequentially CM.
\end{proof} 

Now, our goal is to show that in case of sequentially CM modules it is possible to give formulas for several degree 
functions in terms of certain extension modules. 
At first we compute the homological degree of a module.

\begin{thm}
\label{scm_hdeg}
Let $M$ be a sequentially CM $S$-module of dimension $d$. 
Then we have 
$$
\hdeg M = \deg M  + \sum_{i=0}^{d-1} \binom{d-1}{i} \deg \Ex {n-i} S M .
$$
Moreover, if
$M=F/U$ is a representation of $M$ where  
$F$ is a finitely generated graded free $S$-module and
$U \subset F$ is a graded submodule, then
$$
\hdeg F/U = \hdeg F/\gin(U). 
$$
\end{thm}

\begin{proof}
By the definition of the homological degree we have that 
$$
\hdeg M =\deg M  + \sum_{i=0}^{d-1} \binom{d-1}{i} \hdeg \Ex {n-i} S M.
$$
Since $\hdeg M = \deg M$ for every CM module, the first claim follows. 

If $M=F/U$, then
$$
\hdeg M=
\deg F/U  + \sum_{i=0}^{d-1} \binom{d-1}{i} \deg \Ex {n-i} S {F/U}.
$$
By Theorem 3.1 in \cite{HESB}
the Hilbert functions of the graded modules 
$\Ex {n-i} S {F/U}$ and $\Ex {n-i} S {F/\gin(U)}$
coincide for all $i$. In particular, these modules have the same degree.
Since $F/\gin(U)$ is sequentially CM by \ref{lem-borel-t-is-seq-CM}, this proves the 
second assertion.
\end{proof} 

We will see in Section \ref{sec-bounds} that the statement is not true for an arbitrary $S$-module. 

For a first application of the theorem we need the following result. 

\begin{lem}
\label{help_sec}
Let $M$ be a finitely generated graded $S$-module. 
If $y \in S_1$ 
is $M$-regular and $\Ex {i} S M$-regular for all $i$,
then the module 
$\Ex {i} S M$ is CM of dimension $j$ if and only if $\Ex {i+1} S {M/yM}$ is CM 
of dimension $j-1$.
In this case we have $\deg \Ex {i} S M= \deg \Ex {i+1} S {M/yM}$.
\end{lem}

\begin{proof}
The long exact sequence derived from the short exact sequence
$$
0 \to M(-1) \to M \to M/yM \to 0
$$
splits into short exact sequences 
$$
0 \to \Ex {i} S {M} \to \Ex {i} S {M} (+1)\to \Ex {i+1} S {M/yM} \to 0
$$
from which the assertion follows.
\end{proof}

In \cite[Conjecture 9.4.1]{VA} Vasconcelos conjectured  that
for every $M$-regular element $y\in S$ we have that
$\hdeg M \geq  \hdeg M/yM$. If $M$ is  sequentially CM, then this is true and 
moreover 
we can compute the difference $\hdeg M - \hdeg M/yM$.

\begin{cor}
Let $M$ be a sequentially CM $S$-module of dimension $d$. 
 If $y \in S_1$ is $M$-regular, then
$$
\hdeg M =  \hdeg M/yM  + \sum_{i=1}^{d-2} \binom{d-2}{i} \deg \Ex {n-i} S {M} 
\geq \hdeg M/yM.
$$
In particular,  $\hdeg M =  \hdeg M/yM$ if $\dim M \leq 2$.
\end{cor}

\begin{proof}
It follows from the local duality theorem, that
a prime ideal $P$ of height $i$ is associated to $M$ if and only if $\Ex {i} S 
M_P \neq 0$.
Since $M$ is sequentially CM and therefore $\Ex {i} S M$ is zero or CM of 
dimension
$n-i$,  the associated prime ideals of $\Ex {i} S M$
are exactly the associated prime ideals of $M$ of height $i$.
We deduce that $y$ is also an $\Ex {i} S M$-regular element for all $i$.
The long exact Ext-sequence derived from the short exact sequence
$$
0 \to M(-1) \overset{y}{\to} M \to M/yM \to 0
$$
provides, for all $i<n$,  short exact sequences of the form
$$
0 \to \Ex {n-i} S M \overset{y}{\to} \Ex {n-i} S M (+1) \to \Ex {n-i+1} S {M/yM} 
\to 0. 
$$
 Note that $\Ex {n} S M =0$ because $\depth M>0$. 
 
Now, we show the claim by induction on $d$.  
If $d = 1$ there is nothing  to prove. Let $d\geq 2$.
Observe that $\deg M= \deg M/yM$ and $\dim M/yM=d-1$. 
It follows 
\begin{eqnarray*}
\hdeg M
&=& \deg M     + \sum_{i=1}^{d-1} \binom{d-1}{i} \deg \Ex {n-i} S M \\
&=& \deg M/yM  + \sum_{i=1}^{d-1} \binom{d-1}{i} \deg \Ex {n-i+1} S {M/yM} \\
&=& \deg M/yM  + \sum_{i=1}^{d-1} \biggl( \binom{d-2}{i}+ \binom{d-2}{i-
1}\biggr) \deg \Ex {n-i+1} S {M/yM}\\
&=& \hdeg M/yM  + \sum_{i=1}^{d-1} \binom{d-2}{i} \deg \Ex {n-i+1} S {M/yM}\\
&=& \hdeg M/yM  + \sum_{i=1}^{d-2} \binom{d-2}{i} \deg \Ex {n-i} S {M}\\
\end{eqnarray*}
which is the desired formula.
\end{proof}

Next, we consider the smallest extended degree of a sequentially CM 
module. To this end we recall some well-known results.
For the convenience of the reader
we reproduce the short proofs. 

\begin{lem}
\label{help_cm}
Let $M$ be a finitely generated graded $S$-module and $d=\dim M.$ 
Then $\deg M =\deg \Ex {n-d} S {M}$.
\end{lem}
\begin{proof}
Denote by $P_M$ and $H_M$ the Hilbert polynomial and the Hilbert function of $M$, 
respectively. 
There is the following formula of Serre (cf., e.g., \cite[Theorem 4.4.3]{BRHE98}) 
$$
H_M(j) - P_M(j) 
= \sum_{i=0}^d (-1)^{i} \dim_K  H^i_{\mm}(M)_j 
= \sum_{i=0}^d (-1)^{i} \dim_K  \Ex {n-i} S {M}_{-j}.
$$
Since $\dim \Ex {n-i} S {M} \leq i$ and $\dim \Ex {n-d} S {M} = d$, 
the claim follows by considering $P_M$
and the Hilbert polynomial of $\Ex {n-d} S {M}$ for integers $j \ll 0$. 
\end{proof}

\begin{lem}
\label{help_depth}
Let $M$ be a finitely generated graded $S$-module and $L=H^0_{\mm}(M)$. Then
$$
\deg \Ex {i} S {M/L}=
\begin{cases}
\deg \Ex {i} S {M}  & \text{ for } i<n,\\
0                   & \text{ for } i=n.
\end{cases}
$$
\end{lem}
\begin{proof}
The long exact Ext-sequence derived from the short exact sequence
$$
0 \to L \to M \to M/L \to 0
$$
and the fact that $\Ex {i} S {L}=0$ for $i \neq n$ and $l(\Ex {n} S {L})=l(L)$
imply the assertion. 
\end{proof}

Now, we are ready for the computation of the smallest extended degree. 

\begin{thm}
\label{scm_sdeg}
Let $M$ be a sequentially CM $S$-module of dimension  $d$. Then we have 
$$
\sdeg M 
= \sum_{i=0}^{d} \deg \Ex {n-i} S M = \deg M  + \sum_{i=0}^{d-1} \deg \Ex {n-i} 
S M.
$$
\end{thm} 

\begin{proof} 
Lemma \ref{help_cm} yields the second equality. 
We show the first equality by induction on $d$.
If $d=0$, then $M$ is CM and we have
that $\sdeg M= \deg M=\deg \Ex {n-d} S M$
where the last equality follows from \ref{help_cm}.
Assume that $d>0$.
We  consider two cases. 

(i): Assume $\depth M >0$. 
Since $M$ is sequentially CM 
we can choose an element $y\in S_1$ which 
is $M$-regular and $\Ex {n-i} S M$-regular for all $i$.
It follows 
$$
\sdeg(M) 
= \sdeg(M/yM)  
=   \sum_{i=0}^{d-1} \deg \Ex {n-i} S {M/yM} 
$$
$$
= \sum_{i=0}^{d-1} \deg \Ex {n-i-1} S M   
=  \sum_{i=1}^{d} \deg \Ex {n-i} S M
=  \sum_{i=0}^{d} \deg \Ex {n-i} S M
$$
where the second equality follows from the induction hypothesis and
the third from Lemma \ref{help_sec}.

(ii): Assume $\depth M =0$. 
One of the properties of the smallest extended degree provides 
$$
\sdeg M =\sdeg M/H^0_{\mm}(M)  + l(H^0_{\mm}(M)).
$$
Note that $l(H^0_{\mm}(M))=\deg \Ex {n} S M$ by graded local duality.
Applying case (i) to the module $M / H^0_{\mm}(M)$
and  using Lemma \ref{help_depth} we get 
\begin{eqnarray*}
\sdeg M
&=& \sdeg M/H^0_{\mm}(M) + l(H^0_{\mm}(M))\\
&=& \sum_{i=1}^d \deg \Ex {n-i} S {M/H^0_{\mm}(M)} +  \deg \Ex {n} S M\\
&=& \sum_{i=0}^d \deg \Ex {n-i} S M.\\
\end{eqnarray*} 
This completes the proof. 
\end{proof} 

Finally, we consider the arithmetic degree.

\begin{thm}
\label{scm_adeg}
Let $M$ be a finitely generated graded $S$-module with $d=\dim M$.
If for all $i$ the module $\Ex {n-i} S M$ is zero or CM, then
$$
\adeg M = \sum_{i=0}^{d} \deg \Ex {n-i} S M = \deg M  + \sum_{i=0}^{d-1} \deg 
\Ex {n-i} S M.
$$
In particular, this formula is true for every sequentially CM module.
\end{thm} 

\begin{proof}
By \cite[Proposition 9.1.2]{VA} we know 
 that
$$
\adeg M = \sum_{i=0}^{n} \deg \Ex {i} S {\Ex {i} S M}. 
$$
Since $\Ex {i} S M$ is zero or CM,   
 Lemma \ref{help_cm} provides 
$$
\deg \Ex {i} S {\Ex {i} S M} = \deg \Ex {i} S M.
$$
Using $\Ex {i} S M=0$ for $i<n-d$ and $\deg M= \deg \Ex {n-d} S M$ (by Lemma \ref{help_cm}), we get the claimed equalities.  
In order to conclude the proof, we note that a 
sequentially CM module satisfies the assumption of the theorem.
\end{proof}

Theorem \ref{scm_hdeg}, Theorem \ref{scm_sdeg} and Theorem \ref{scm_adeg}
imply in particular the following result.

\begin{cor}
\label{sec_cor2}
Let $M$ be a finitely generated graded $S$-module which is sequentially CM. 
Then we have
$$
\deg M \leq \adeg M = \sdeg M \leq \hdeg M.
$$
Furthermore, 
\begin{enumerate}
\item
$\deg M =\adeg M$ if and only if $M$ is Cohen-Macaulay.
\item
$\sdeg M = \hdeg M$ if and only if  
$\deg \Ex {n-i} S M =0$ for $i=1,\dots,d-2$.
\end{enumerate}
\end{cor}

In Section \ref{sec-bounds} we will see that some of these relations are true in much greater generality.

%
%
%
\section{Buchsbaum modules} \label{sec_buchsbaum} 

Recall that 
$S=K[x_1,\dots,x_n]$ is the polynomial ring over the infinite field $K$ with 
graded maximal ideal $\mm=(x_1,\dots,x_n)$.
Let $M$ be a finitely generated graded $S$-module of dimension $d=\dim M$.
The module $M$ is called a {\em Buchsbaum module}
if $d=0$ or $d>0$ and every homogeneous system of parameters $y_1,\dots,y_d$ of 
$M$
is a {\em weak $M$-sequence}, i.e.
$$
(y_1,\dots,y_{i-1})M\colon y_i =(y_1,\dots,y_{i-1}) M \colon \mm \text{ for } 
i=1,\dots,d.
$$ 
We need some properties of Buchsbaum modules.
If $M$ is a Buchsbaum module, then:
\begin{enumerate}
\item
$\mm H_\mm^i(M)=0$ for $i=1,\dots,d-1$;
\item
If $y_1,\dots,y_r$ is part of a homogeneous system of parameters of $M$,
then also $M/(y_1,\dots,y_r)M$ is a Buchsbaum module of dimension $d-r$.
\end{enumerate}
Note that property (i) implies that 
$l(H_\mm^i(M))<\infty$. Thus, local duality provides that 
$l(H_\mm^i(M))=l(\Ex {n-i} S M)$
and $\mm \Ex {n-i} S M=0$.
An arbitrary module $M$ that only satisfies property (i) is called
a {\em quasi-Buchsbaum module}. For more details on the theory of Buchsbaum 
modules
we refer to the book of St\"uckrad and Vogel \cite{STVO}.

In this section we study the behavior of degree functions when they are applied to Buchsbaum 
modules. In case of the homological degree, the 
 following result was already noted in \cite[Theorem 9.4.1]{VA}. It follows immediately from the definition of $\hdeg$ 
because every module $M$ of finite length satisfies $\hdeg M = l(M)$. 

\begin{prop}
\label{buchs_hdeg}
Let $M$ be a finitely generated graded $S$-module, $d=\dim M$
such that $l(\Ex {n-i} S M)<\infty$ for $i=0,\dots,d-1$.
Then 
$$
\hdeg M 
= 
\deg M  + \sum_{i=0}^{d-1} \binom{d-1}{i} l(\Ex {n-i} S M).
$$
In particular, this formula is true for every Buchsbaum module.
\end{prop}


Next we compute the smallest extended degree of a Buchsbaum module. The 
 theorem below 
was first stated in Gunston's thesis \cite[Proposition 3.2.3]{GU}, but with the weaker hypothesis that $M$ is quasi-Buchsbaum. 
However Gunston's  proof does not work in this generality, because if $M$ is a quasi-Buchsbaum module
and $y$ is a homogeneous parameter element for $M$, then  $M/yM$ is in general not a quasi-Buchsbaum module.
For an example of such a module see \cite[Example 7.4]{NA}. 
Since Gunston's result is not published elsewhere, we give a proof of the theorem.

\begin{thm}
\label{buchs_sdeg}
Let $M$ be a  graded Buchsbaum $S$-module of dimension $d$.
Then we have  
$$
\sdeg M 
= 
\deg M  + \sum_{i=0}^{d-1} \binom{d-1}{i} l(\Ex {n-i} S M)
=
\hdeg M.
$$
\end{thm}

\begin{proof}
We use induction on $d$.
If $d=0$, then $M$ is CM and we have that $\sdeg M =\hdeg M =\deg M$. 
Let $d>0$. Assume that 
 $\depth M =0$ and that the assertion is already shown for modules with positive 
depth.  Then we get because $M/H_\mm^0(M)$ is Buchsbaum, too, that 
\begin{eqnarray*}
\sdeg M &=& 
\sdeg M/H_\mm^0(M) + l(H_\mm^0(M))\\
&=&  
\deg M/H_\mm^0(M)  + \sum_{i=1}^{d-1} \binom{d-1}{i} l(\Ex {n-i} S 
{M/H_\mm^0(M)}) + l(\Ex {n} S {M})\\
&=&  
\deg M  + \sum_{i=0}^{d-1} \binom{d-1}{i} l(\Ex {n-i} S {M}).
\end{eqnarray*}
Here the first equality follows from the properties of $\sdeg$, the second equality 
since 
$\depth M/H_\mm^0(M)>0$ (note that $l(\Ex {n} S {M/H_\mm^0(M)})=0$),
and the third one from Lemma \ref{help_depth}. 

It remains to consider the case $\depth M >0$. Choose an $M$-regular element $y\in 
S_1$ and consider the short exact sequence
$$
0 \to M(-1) \overset{y}{\to} M \to M/yM \to 0.
$$ 
Observe that $y\cdot \Ex {n-i} S {M}=0$. Hence the associated long exact $\Ext$-sequence
splits into short exact sequences of the form
$$
0 \to \Ex {n-i-1} S {M}(+1) \to \Ex {n-i} S {M/yM} \to \Ex {n-i} S {M} \to 0.
$$
Thus, we get that $l(\Ex {n-i} S {M/yM})= l(\Ex {n-i-1} S {M}) + l(\Ex {n-i} S {M})$.
Since  $M/yM$ is again a Buchsbaum module we may apply the induction hypothesis to 
it and obtain 
\begin{eqnarray*}
\sdeg M &=& 
\sdeg M/yM\\
&=&  
\deg M/yM  + \sum_{i=0}^{d-2} \binom{d-2}{i} l(\Ex {n-i} S {M/yM})\\
&=&  
\deg M  + \sum_{i=0}^{d-2} \binom{d-2}{i} \biggl(l(\Ex {n-i-1} S {M}) + l(\Ex 
{n-i} S {M})\biggr)\\
&=&  
\deg M  + \sum_{i=0}^{d-1} \binom{d-1}{i} l(\Ex {n-i} S {M}).
\end{eqnarray*}
Comparing with \ref{buchs_hdeg} we see that $\sdeg M = \hdeg M$ and this concludes the 
proof.
\end{proof}

There is also a formula in case of the arithmetic degree. 

\begin{prop}
\label{buchs_adeg}
Let $M$ be a  graded Buchsbaum $S$-module of dimension $d$.
Then we have
$$
\adeg M = \sum_{i=0}^{d} l(\Ex {n-i} S M) = \deg M  + \sum_{i=0}^{d-1} l(\Ex {n-
i} S M).
$$
\end{prop}
\begin{proof}
This follows from \ref{scm_adeg} since all modules $\Ex {n-i} S M$ are zero or 
CM of dimension $0$.
\end{proof}

Combining the previous results we get a statement that is similar to \ref{sec_cor2}. 

\begin{cor}
\label{buchs_cor1}
Let $M$ be a graded Buchsbaum $S$-module of dimension  $d$.
Then we have
$$
\deg M \leq \adeg M \leq \sdeg M = \hdeg M.
$$
Furthermore,  
\begin{enumerate}
\item
$\deg M =\adeg M$ 
if and only if 
$M$ is Cohen-Macaulay.
\item
$\adeg M = \sdeg M$ 
if and only if 
$\deg \Ex {n-i} S M =0$ for $i=1,\dots,d-2$.
\end{enumerate}
\end{cor}

Buchsbaum modules form another class of modules where we can  give an affirmative answer
to Conjecture 9.4.1 in \cite{VA}. 

\begin{cor}
Let $M$ be a graded Buchsbaum $S$-module. 
If  $y \in S_1$ is $M$-regular, then
$$
\hdeg M =  \hdeg M/yM.
$$
\end{cor}
\begin{proof}
We have that $\hdeg M = \sdeg M = \sdeg M/yM = \hdeg M/yM$
where we used Corollary \ref{buchs_cor1} and the fact that $M/yM$ is again a 
Buchsbaum module.
\end{proof}

%
%
%
\section{Bounds for degree functions} \label{sec-bounds} 

We apply the results of the last sections to compare degree functions and to study their behavior when passing to certain 
monomial modules. This leads to various bounds. 

The starting point is the following refinement of \cite[Proposition 9.4.2]{VA}. 

\begin{thm}
\label{sec_cor3}
Let $M$ be a finitely generated graded $S$-module. 
Then we have 
$$
\deg M \leq \adeg M \leq \sdeg M \leq \hdeg M.
$$
\end{thm} 

Using our previous results we can give a new, more conceptual proof. We need an extension of Theorem 3.2 in \cite{STTRVO} to modules. 


\begin{lem} \label{lem-adeg-gin} 
If $U$ is a finitely generated graded submodule of the free graded $S$-module $F$ then we have 
$$
\adeg F/U \leq \adeg F/\gin (U). 
$$
\end{lem} 

\begin{proof} 
We use again the formula 
$$
\adeg F/U = \sum_{i=0}^{n} \deg \Ex {i} S {\Ex {i} S {F/U}}. 
$$
Observe that $\Ex {i} S {\Ex {i} S {F/U}} \neq 0$ if and only if $\dim {\Ex {i} S {F/U}} = n-i$. Thus, we get
$$
\deg \Ex {i} S {\Ex {i} S {F/U}} \leq \deg {\Ex {i} S {F/U}} \leq \deg {\Ex {i} S {F/\gin (U)}}
$$
where the first estimate is a consequence of \ref{help_cm} and the second inequality follows from the fact that, by \cite{SB}, 
$\dim_K \Ex {i} S {F/U}_j \leq \dim_K \Ex {i} S {F/\gin(U)}_j$ for all integers $i, j $. 
\end{proof} 

\begin{proof}[Proof of Theorem \ref{sec_cor3}]
We know that $\sdeg M \leq \hdeg M$ is true by the properties of the smallest 
extended degree.
Choose a presentation $M=F/U$
where $F$ is a finitely generated free graded $S$-module
and $U$ is a graded submodule of $F$.
Then we get
$$
\adeg M= \adeg F/U \leq \adeg F/\gin(U) = \sdeg F/\gin(U) $$
$$= \sdeg F/U = \sdeg M \leq \hdeg M.
$$
Here the second inequality follows from Lemma \ref{lem-adeg-gin}
and the third one from \ref{sec_cor2} because $F/\gin(U)$ is sequentially CM. 
\end{proof}

In order to estimate $\sdeg M$ and $\hdeg M$ it seems natural by now to consider a presentation $M = F/U$ where $F$ is a 
free module and to compare the degrees of $M$ with the ones of $F/\gin (U)$ and $F/U^{lex}$. This works well to give a lower bound.

\begin{prop}
\label{hdeg_lower}
Let $M$ be a finitely generated graded $S$-module. 
Let $M=F/U$ be a representation of $M$ where  
$F$ is a finitely generated graded free $S$-module and
$U \subset F$ is graded submodule.  
Then we have 
$$
\hdeg M \geq \sdeg M= \deg M + \sum_{i=0}^{d-1} \deg \Ex {n-i} S {F/\gin(U)}
$$
and equality is true if $M$ is a Buchsbaum module.
\end{prop}
\begin{proof}
We have that
$$
\hdeg M 
\geq \sdeg M
= \sdeg F/U
= \sdeg F/\gin(U)
= \sum_{i=0}^d \deg \Ex {n-i} S {F/\gin(U)}.
$$
Here the inequalities and equalities follow from the
properties of $\sdeg$ (see the remarks after \ref{sdeg_define})
and Theorem \ref{scm_sdeg} because $F/\gin(U)$ is sequentially CM by 
\ref{boreltype}. 
If $M$ is a Buchsbaum module, then Corollary \ref{buchs_cor1} shows the claimed equality.
\end{proof}

In order to give an upper bound for $\hdeg M$  we have to restrict ourselves to certain classes of modules 
because we show in Section  \ref{sec-examples} that the analogous result is not true for an arbitrary graded $S$-module. 


\begin{thm}
\label{lex_bound}
Let $M$ be a finitely generated graded $S$-module
which is sequentially CM or a Buchsbaum module. 
Let $M=F/U$ be a representation of $M$ where  
$F$ is a finitely generated graded free $S$-module and
$U \subset F$ is graded submodule.  
Then we have
$$
\hdeg M  \leq  \hdeg F/U^{lex}. 
$$
\end{thm} 

\begin{proof}
Consider first the case that $M=F/U$ is sequentially CM.
Sbarra proved in his thesis (see \cite{SB} for a published proof) 
that 
$$
\dim_K \Ex {n-i} S {F/U}_j \leq \dim_K \Ex {n-i} S {F/U^{lex}}_j
\text{ for all } i,j.
$$
Since $U^{lex}$ is of Borel-type, it is  sequentially CM by Lemma \ref{lem-borel-t-is-seq-CM} .
Thus, the modules 
$\Ex {n-i} S {F/U}$ and $\Ex {n-i} S {F/U^{lex}}$
are zero or CM of dimension $i$.
Using the inequalities above, this implies that 
$$
\deg\Ex {n-i} S {F/U} \leq \deg \Ex {n-i} S {F/U^{lex}}
\text{ for all } i. 
$$ 
Now Theorem \ref{scm_hdeg} shows  that
$\hdeg F/U  \leq  \hdeg F/U^{lex}.$

Second, assume that $F/U$ is Buchsbaum. 
Then we know from \ref{buchs_cor1} that
$\hdeg F/U = \sdeg F/U$.
Recall that  $\sdeg F/U= \sdeg F/\gin(U)$. Applying 
Theorem \ref{scm_sdeg} to
$F/\gin(U)$ and $F/\gin(U)^{lex}=F/U^{lex}$ and using an argument analogous to the one above in the case of $\hdeg$ of sequentially CM modules, we obtain 
$$
\sdeg F/\gin(U) \leq \sdeg F/\gin(U)^{lex}= \sdeg F/U^{lex}.
$$ 
It follows that 
$$
\hdeg F/U \leq \sdeg F/U^{lex} \leq \hdeg F/U^{lex}.
$$
This concludes the proof.
\end{proof}

Bounding
$\sdeg$ is much easier as the proof of
Theorem \ref{lex_bound} shows.
\begin{thm}
\label{sdeg_lex}
Let $M$ be a finitely generated graded $S$-module. 
Let $M=F/U$ be a representation of $M$ where  
$F$ is a finitely generated graded free $S$-module and
$U \subset F$ is graded submodule.  
Then
$$
\sdeg M = \sdeg F/\gin(U) \leq \sdeg F/ U^{lex}.
$$
\end{thm}
\begin{proof}
We know already the first equality.
Now the second part of the proof of
\ref{lex_bound} shows for an arbitrary graded submodule $U \subset F$ that
$$
\sdeg F/\gin(U) \leq \sdeg F/ \gin(U)^{lex}=\sdeg F/ U^{lex}.
$$
\end{proof}

Note that formulas for the bounds for $\hdeg F/U^{lex}$ and $\sdeg F/U^{lex}$ are given by \ref{scm_hdeg} and \ref{scm_sdeg}. 
Thus, getting effective estimates amounts to computing degrees of certain extension modules. This can be done efficiently. 

Indeed, observe that $\gin(U)$ and $U^{lex}$ are monomial submodules  of Borel-type 
(cf.\ \ref{lem-borel-t-is-seq-CM}).
In Section \ref{computational} we will show that one can fastly compute the degree of
$\Ex {n-i} S {F/V}$ where $V \subset F$ is monomial of Borel-type, 
if one just knows the unique minimal system of monomial generators of $V$. Thus, it is possible to compute our bounds using computer 
algebra systems
like CoCoA \cite{CO}, Macaulay 2 \cite{MC2} or Singular \cite{Sing}. 

%
%
%
\section{Counterexamples }  \label{sec-examples} 

The work on this paper started by trying to prove the
following conjecture (see Gunston \cite[Conjecture 2.5.3]{GU} and the book of 
Vasconcelos \cite[page 262]{VA}): One of the relations 
$$
\hdeg S/I \leq \hdeg S/\gin(I) \text{ or } \hdeg S/I \geq \hdeg S/\gin(I)  
$$ 
is true for {\em every} homogeneous ideal $I \subset S$. 

Very often the lexicographic ideal $I^{lex}$ associated to $I$ has extremal 
properties 
with respect to invariants of $I$. Therefore it is natural to study the related 
problem: Is one of the relations 
$$
\hdeg S/I \leq \hdeg S/I^{lex} \text{ or } \hdeg S/I \geq \hdeg S/I^{lex} 
$$ 
true for all ideals $I$? We have seen that the analogous problem for $\sdeg$ has 
a positive answer because in  Theorem \ref{sdeg_lex} we proved that 
$$
\sdeg S/I = \sdeg S/\gin(I) \leq \sdeg S/ I^{lex}. 
$$
Now, we will show that for $\hdeg$ all inequalities are false in general. 
First, we consider the comparison of $I$ and $I^{lex}$. 

\begin{ex} 
\label{sec_ex}
Let $d\geq 3$, $g < \binom{d-2}{2}$ be integers and set  $a := \binom{d-1}{2}-g$.  Consider the ideal  
$$
I=(x^2 , xy , y^d , y^{d-1}z^{a-d+2} + xt^{a}) \subset K[x,y,z,t] =: S.
$$
The ideal $I$ is the homogeneous ideal of an extremal projective curve of degree $d$ and genus $g$
as considered in \cite[Example 4.5]{NAhart} 
(our $a$ is the one of that Example plus $d-2$!). In \cite{NAhart}, it is shown that 
$$
H^1_\mm(S/I) \cong K[x,y,z,t]/(x,y,z^{a-d+2}, t^a) (d-1-a), 
$$
thus 
$$
\deg \Ex {3} S {S/I}=l(H^1_\mm(S/I))=a(a-d+2)
$$
Since $\dim S/I =2$, $\depth S/I >0$ 
(for example, the element $t$ is $S/I$-regular), we get
$$
\hdeg S/I = d + a(a-d+2).
$$

Next, we compute $\hdeg S/I^{lex}$. The saturation of $I^{lex}$ is (cf., e.g., \cite{Bayer}) 
$$
(I^{lex})^{sat}=(x,y^{d+1},y^dz^a)
$$
Using, e.g., \cite[Corollary 2.6]{HEPOVL} or \cite[Lemma 3.4]{NA03} we obtain 
$$
\deg \Ex {3} S {S/I^{lex}}=\deg \Ex {3} S {S/(I^{lex})^{sat}}= a.
$$
By considering the Hilbert function of $I$ 
and using again the theorem of Serre (\cite[Theorem 4.4.3]{BRHE98}))
one gets
$$
\deg \Ex {4} S {S/I^{lex}}=l(H^0_\mm(S/I^{lex}))=a(d-1).
$$
It follows that 
$$
\hdeg S/I^{lex}=d + a + a(d-1)=d+ad
$$
For $d=3, a=2$ we obtain
$\hdeg S/I= 5 < 9=\hdeg S/I^{lex}$
and for $d=3, a=5$ we get 
$\hdeg S/I= 23 > 18=\hdeg S/I^{lex}$.
Therefore, in general, there is no relation  between
$\hdeg S/I$ and $\hdeg S/I^{lex}$, not even for two-dimensional rings. 
\end{ex} 

Now, we turn to the comparison of $I$ and $\gin (I)$. 

\begin{ex} \label{ex-c-gin} 
First, we consider again the ideal $I$ of \ref{sec_ex}. Its generic initial ideal has been computed in \cite[Proposition 5.5]{NAhart}. 
It is 
$$
\gin(I)=(x^2,xy,y^d,y^{d-1}z^{a-d+2}). 
$$
Using \cite[Corollary 2.6]{HEPOVL} we get 
$$
\deg \Ex {3} S {S/\gin(I)}=l(H^1_\mm(S/I))= a-d+2.
$$
Thus, since $\depth S/\gin (I) = 1$ we obtain 
$$
\hdeg S/ \gin(I) = d + a-d+2=a+2.
$$
Hence, for $d=3$ and $a=2$ we get $\hdeg S/I= 5    > 4=\hdeg S/\gin(I)$.

Second, we take a graded Buchsbaum ring $S/J$ of dimension $d\geq 3$
over a polynomial ring $K[x_1,\dots,x_n]$ with $\Ex {n-1} S {S/ J} \neq 0 $. 
Then we have 
$$
\hdeg S/J = \sdeg S/J = \sdeg S/\gin(J) < \hdeg S/\gin(J).
$$
The first equality follows from Theorem \ref{buchs_cor1} and
the second is a property of $\sdeg$.
The third inequality is a consequence of \ref{sec_cor2}, the fact
that $S/\gin(J)$ is sequentially CM by 
\ref{lem-borel-t-is-seq-CM},  and 
that $\Ex {n-1} S {S/\gin(J)} \neq 0 $.
The latter we deduce from  Sbarra's result  
in \cite{SB},  
$
\dim_K \Ex {i} S {S/J} \leq  \Ex {i} S {S/\gin(J)}_j \text{ for all } i,j.
$

This shows that, in general, there is no relation  between
$\hdeg S/I$ and $\hdeg S/\gin (I)$. 
\end{ex}

%
%
%
\section{Algorithms}
\label{computational}

We have seen that the smallest extended degree has a number of 
 nice properties that are not shared by the homological degree. However, the homological degree has the advantage 
that is defined by an explicit formula. The goal of this section is to   
 present an
algorithm which shows that it is possible to compute  effectively the smallest extended degree 
by using computer algebra systems 
like CoCoA \cite{CO}, Macaulay 2 \cite{MC2}, or Singular \cite{Sing}.

The idea for the computation of the smallest extended degree is to use the fact
that
$\sdeg F/U = \sdeg F/\gin(U)$ for a graded submodule $U$ 
of a finitely generated graded free $S$-module $F$. This relies on the efficient computation of 
$\deg \Ex {i} S {F/U}$ whenever $U$ is of Borel-type. 
The key result is the following algorithm
(see Sturmfels, Trung and Vogel \cite{STTRVO} for a related result concerning
$\adeg S/I$ where $I$ is a monomial ideal):

\begin{alg}
\label{alg_sdeg2}
Let $U$ be a monomial submodule of a finitely generated graded free $S$-module 
$F$
with homogeneous basis $e_1,\dots,e_m$.
Assume that $U$ is a of Borel-type, i.e.
$$
U: x_i^\infty = U: (x_1,\dots,x_i)^\infty \text{ for } i=1,\dots,n. 
$$
Define inductively graded submodules $U_0,\dots,U_n$ of $F$ 
as follows:
\begin{enumerate}
\item
Set $U_0=M$.
\item
Assume that
$U_0,\dots,U_{i-1}$ are chosen. 
Put $U_i=U_{i-1}\colon x_{n-i+1}^\infty.$ 
\end{enumerate}
Let $G_i$ be the unique minimal system of monomial generators of the monomial 
module $U_i$.
Then for all $x^ue_j \in G_i$ we have that $x_t \nmid x^u$ for $t \geq n-i+1$.
Let $V_i$ be the monomial submodule of $F$ generated by $G_i$ as a 
$K[x_1,\dots,x_{n-i}]$-module
and denote by $V_i^{sat}$ the saturation of this module as a 
$K[x_1,\dots,x_{n-i}]$-module.
Then $V_i^{sat}/V_{i}$ is of finite length for $i=0,\dots,n$
and we have 
$$
\deg \Ex {n - i} S {F/U} = l((V_i^{sat}/V_{i}), 
$$
thus
$$
\sdeg F/U =  \sum_{i=0}^{n} l(V_i^{sat}/V_{i}).
$$
\end{alg}
\begin{proof}
It follows from \ref{scm_sdeg} that
$$
\sdeg F/U = \sum_{i=0}^{n} \deg \Ex {i} S {F/U}
$$
because $F/U$ is sequentially CM by \ref{lem-borel-t-is-seq-CM}
and $\Ex {i} S {F/U}=0$ for $i<n-\dim F/U$.
We have that $U=\Dirsum_{i=1}^m I_je_j$ for monomial ideals $I_j \subset S$ of 
Borel-type.
All computations in \ref{alg_sdeg2} and the Ext-modules commute with
finite direct sums. 
Therefore we may assume that $F/U=S/I$ for some monomial ideal $I \subset S$ of 
Borel-type.
Now the claim follows from the structure theorem of $\Ex {i} S {S/I}$ in 
\cite[Corollary 2.6]{HEPOVL}.
Note that in the proofs for Corollary 2.6 in \cite{HEPOVL}, Herzog, Popescu and 
Vladoiu
consider the unique shortest chain of monomial ideals. We changed their computations 
a little bit.
In our notation it can happen that $U_i=U_{i-1}$, but then $V_{i}^{sat}=V_{i}$ 
and $l(V_i^{sat}/V_{i})=0$.
\end{proof}

\begin{rem}
Algorithm \ref{alg_sdeg2} is very fast.
If we have computed the unique minimal systems $G_i$ of monomial generators 
of $U_i$,
then it is easy to determine 
$$
U_i=U_{i-1}\colon x_{n-i+1}^\infty, 
$$ 
because we have
$$
G_{i+1}=\{ x^u/x_{n-i+1}^{m_{n-i+1}(u)}  e_j\colon x^ue_j \in G_i \}
$$
where for a monomial $x_u$ we set $m_t(u)= \max \{s\colon x_t^s|x^u \}$.
The computation of the length of a quotient of monomial submodules can also  very 
efficiently be done. 
\end{rem}

Now one can try to extend \ref{alg_sdeg2} to arbitrary modules as follows.
\begin{alg}
\label{alg_sdeg3} \label{alg-sdeg} 
Let $M$ be a finitely generated graded $S$-module. 
\begin{enumerate}
\item
Determine a presentation 
$M=F/U$ 
where $F$ is a finitely generated graded free $S$-module and
$U$ is a graded submodule of $F$.
\item
Compute $\gin(U)$ with respect to the reverse lexicographic order.
\item
Find $s := \sdeg F/ \gin(U)$ using \ref{alg_sdeg2}.
\end{enumerate}
Then we have that
$\sdeg M =  s.$
\end{alg}
\begin{proof}
This follows from the fact that $\sdeg F/U = \sdeg F/\gin(U)$
and  $\gin(U)$ is a monomial submodule of $F$ of Borel-type.
\end{proof}

\begin{rem}
There is a serious (theoretical) problem with Algorithm \ref{alg_sdeg3}.
To compute $\gin(U)$ one 
usually takes randomly chosen coordinates $y_1,\dots,y_n$ of $S$, 
applies the automorphism $\phi$ of $F$ induced by $x_i \mapsto y_i$  on $U$, 
and then computes $\gin(U)$ as the initial module of $\phi(U)$ 
with respect to the reverse lexicographic term order on $F$.
Since for a ``generic'' $\phi$ we have indeed $\gin (U) = \ini \phi(U)$ this procedure will almost 
always correctly determine $\gin (U)$.  However, there is no criterion to decide if the monomial module one
gets by these computations is in fact the generic initial module
of $U$.  Hence, one cannot be certain if the result is correct.
For practical purposes, the above procedure  is of course good enough, since 
the probability not to get $\gin(U)$ is zero.
\end{rem} 

In his thesis Gunston \cite{GU} has proposed an algorithm for computing $\sdeg M$ that uses general hyperplane sections. 
His procedure also has the theoretical problem that there is no criterium to check if a randomly chosen linear form is indeed 
general enough. If the module $M$ has dimension $d$ then Gunston's algorithm requires the computation of $d+1$ Gr\"obner 
bases whereas our algorithm has the advantage that it needs just one Gr\"obner basis computation.

\end{document}